\documentclass[11pt]{amsart}
\usepackage{amsmath, amssymb, amscd}
\usepackage{epsfig, psfrag, setspace}

\newcommand{\la}{\langle}
\newcommand{\ra}{\rangle}
\newcommand{\Ln}{{{\mathcal L}{\it ie\/}(n)}}
\newcommand{\Lie}{{{\mathcal L}{\it ie\/}}}
\newcommand{\Dn}{{{\mathcal E}{\it il\/}(n)}}
\newcommand{\Eil}{{{\mathcal E}{\it il\/}}}
\newcommand{\Pn}{{{\mathcal P}{\it ois\/}(n)}}
\newcommand{\Op}{{\mathcal O}}  
\newcommand{\mC}{{\mathcal{C}}}
\newcommand{\one}{{\bf{1}}}
\newcommand{\poi}{{{\mathcal P}{\it ois}}}
\newcommand{\DPn}{{{\mathcal S}{\it iop\/}(n)}}
\newcommand{\DP}{{{\mathcal S}{\it iop\/}}}
\newcommand{\Pnk}{{{\mathcal P}{\it ois\/}_k(n)}}
\newcommand{\DPnk}{{{\mathcal S}{\it iop\/}_k(n)}}
\newcommand{\PnP}{{{\mathcal P}{\it ois\/}_P(n)}}
\newcommand{\DPnP}{{{\mathcal S}{\it iop\/}_P(n)}}
\newcommand{\olie}{{{\mathcal L}{\it ie\/}^o}}
\newcommand{\oeil}{{{\mathcal E}{\it il\/}}^o}
\newcommand{\oPn}{{{\mathcal P}{\it ois\/}^o(n)}}
\newcommand{\oDPn}{{{\mathcal S}{\it iop\/}^o(n)}}
\newcommand{\oPnk}{{{\mathcal P}{\it ois\/}^o_k(n)}}
\newcommand{\oDPnk}{{{\mathcal S}{\it iop\/}^o_k(n)}}
\newcommand{\oPnP}{{{\mathcal P}{\it ois\/}^o_P(n)}}
\newcommand{\oDPnP}{{{\mathcal S}{\it iop\/}^o_P(n)}}
\newcommand{\Fo}{{\text{Set}^{op}}}
\newcommand{\mB}{{\mathcal B}}
\newcommand{\n}{\mathbf{n}}

\newtheorem{theorem}{Theorem}[section]
\newtheorem{corollary}[theorem]{Corollary}
\newtheorem{lemma}[theorem]{Lemma}
\newtheorem{proposition}[theorem]{Proposition}

\newtheorem{mydiagram}{Figure}

\theoremstyle{definition}                          
\newtheorem{definition}[theorem]{Definition}  

\newcommand{\refT}[1]{Theorem~\ref{T:#1}}
\newcommand{\refC}[1]{Corollary~\ref{C:#1}}
\newcommand{\refP}[1]{Proposition~\ref{P:#1}}
\newcommand{\refD}[1]{Definition~\ref{D:#1}}
\newcommand{\refL}[1]{Lemma~\ref{L:#1}}

\newcommand{\refF}[1]{Figure~\ref{F:#1}}

\begin{document}

%\doublespacing

\title{A pairing between graphs and trees}
\author{Dev P. Sinha}
%\dedicatory{Dedicated to the memory of my father, Om P. Sinha}
\address{Department of Mathematics, University of Oregon, Eugene OR 97403, USA}

\thanks{This work was partially supported by the National Science Foundation}
\maketitle

In this paper  we develop a canonical pairing between trees
and graphs, which passes to their quotients by
Jacobi and Arnold identities. Our first
main result is that on these quotients the pairing is perfect, which makes
it an effective and simple tool for understanding
the Lie and Poisson operads, providing canonical duals.  Passing from
the operads to free algebras over them, we
get canonical models for cofree Lie coalgebras.
The functionals on free Lie algebras
which result are defined without reference to the embedding
of free Lie algebras in tensor algebras.
In the course of establishing our main results we reprove standard facts
about the modules $\Lie(n)$.  We apply the pairing
to develop product, coproduct and (co)operad structures,
defining notions such as a partition of forests which may be useful elsewhere.  
Remarkably, we find the cooperad which  dual to the Poisson operad more manageable
than the Poisson operad itself.

This pairing arises as the pairing between canonical bases for homology and
cohomology of configurations in Euclidean space.  We  elaborate
on this topology in the expository paper \cite{Sinh06}.
A variant of this pairing first appears in
work of Melancon and Reutenaur on odd-graded free Lie algebras \cite{MeRe96}; 
see Section~\ref{S:freelie}. The pairing was independently developed and applied by 
Tourtchine \cite{Tour98, Tour04};
see further commentary at the end of Section~\ref{S:basic}.  We give a unified, explicit, and
fully self-contained account here to be built on in a number of different
directions in future work, which will include fundamental new results on Lie coalgebras
in algebra and topology \cite{SiWa06}.

%We would like to thank Victor Tourtchine and Christophe Reutenaur for bringing their
%work to our attention.

\section{The (even) Lie configuration pairing}\label{S:basic}

\begin{definition}\label{D:firstdef}
\begin{enumerate}
\item A Tree is an isotopy class of acyclic graph whose vertices
are either trivalent or univalent, with a 
distinguished univalent vertex called the root, embedded in the 
upper half-plane with the root at the origin.  
Univalent non-root vertices are called leaves, and they
are labeled by some set $L$, usually $\n = \{ 1, \ldots, n \}$.
Trivalent vertices are also called internal vertices.
\item The height of a vertex in a Tree is the number of edges between that vertex
and  the root.
\item Define the nadir of a path in a Tree to be the vertex of lowest
height which it traverses.  
\item A Graph is a connected oriented graph with vertices labeled by  
some set $L$, taken to be the appropriate $\n$ unless otherwise noted.  
\item Given a Tree $T$ and a Graph $G$ labeled by the same set, define 
$$\beta_{G,T} : \{ {\rm edges \; of }\; G\} \to \{ {\rm internal\; vertices \; of  }\; T\}$$
by sending an edge $e$ connecting vertices labeled by $i$ and $j$ to the nadir of
the shortest path $p_T(e)$ between the leaves of $T$ labeled $i$ and $j$, which
we call an edge path.
Let $\la G, T \ra_{(2)}$ to
be one if $\beta_{G,T}$ is a bijection and zero otherwise.
\item    In the definition of $\beta_{G,T}$, let $\tau_{G,T} = (-1)^{N}$ where $N$ is the number of
edges $e$ in $G$ for which $p_T(e)$  travels from left to right (according to the half-planar embedding) 
at its nadir.
Define the configuration pairing $\la G, T \ra$ as $\tau_{G,T} \la G, T \ra_{(2)}$.
\end{enumerate}
\end{definition}

See  Figure~\ref{F:basic} for an illustration.  Note that $\la G, T\ra_{(2)}$ 
is defined without reference to the orientation data of the Graph $G$ or the planar
embedding of the Tree $T$.  We may alternately view a Tree through its  set
of vertices ordered by $v \leq w$ if $v$ is in the shortest path between $w$
and the root.  In this language, $\beta_{G,T}$ sends an edge in $G$ to the greatest
lower bound of the two leaves in $T$ labeled by the endpoints of the edge.

%In a Tree each
%internal vertex has one edge closest to the root, which we call its base, and
%two other edges we call branches.   The planar embedding determines and is determined
%by a left-right ordering of the branches of each vertex.  Another way to define the sign
%$\tau$ is that $N$ is the number of edges $e$ such that $p_T(e)$ first crosses the right
%branch of its nadir and then the left.  

\begin{center}
\begin{minipage}{12cm}
{\psfrag{sig}{$\sigma$}
$$\includegraphics[width=12cm]{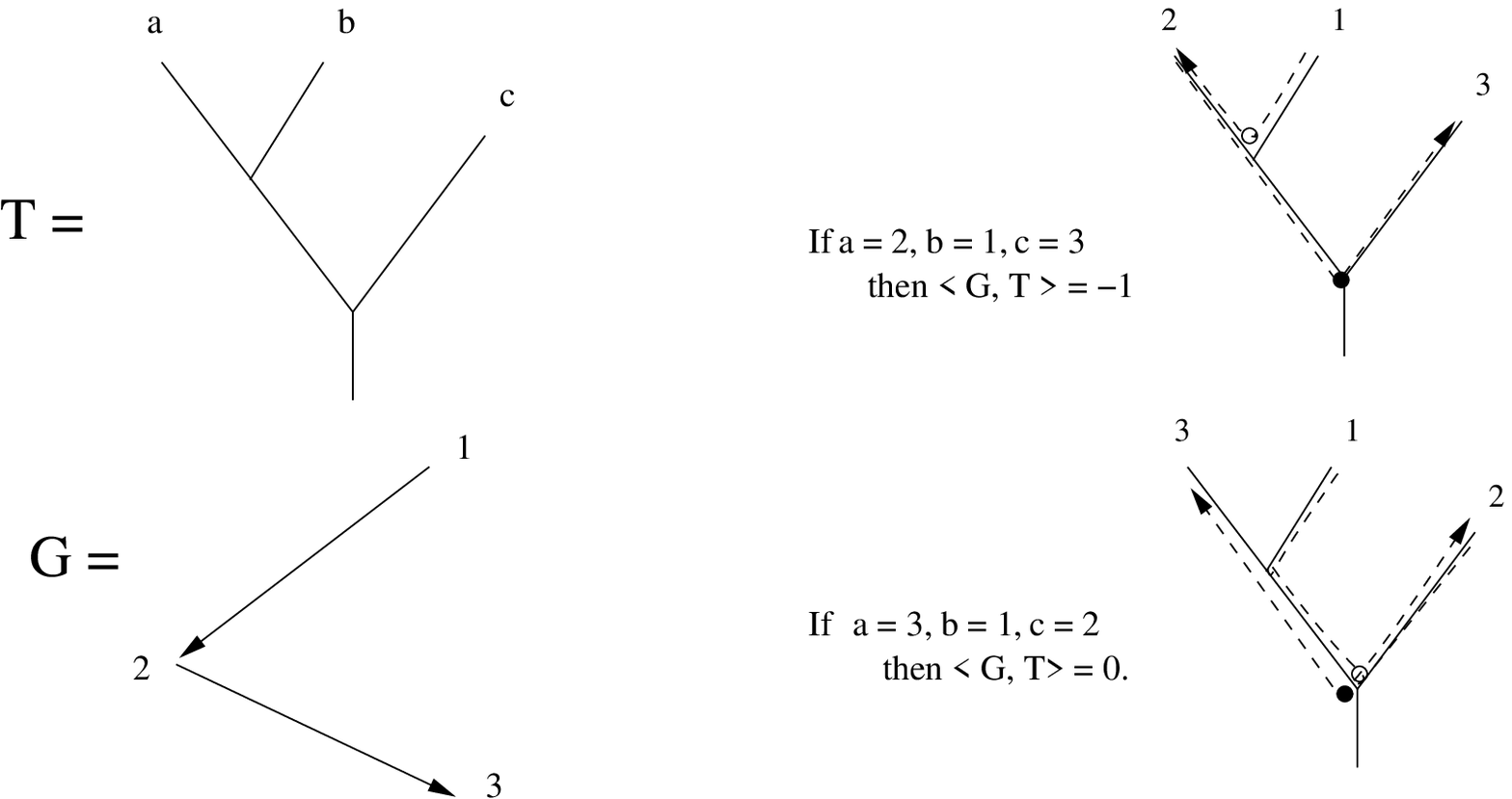}$$}
\begin{mydiagram}\label{F:basic}
Two examples of the configuration pairing, involving one underlying tree
but two different labelings of its leaves.  
\end{mydiagram}
\end{minipage}
\end{center}
\medskip

%\begin{figure}\label{F:basic}
%\end{figure}

The pairing may be 
defined for non-trivalent trees, but we have yet to find an application of such 
generality.    We extend the pairing to free modules over a fixed ground ring
generated by Trees and Graphs, setting notation as follows.

\begin{definition}
Fix a ground ring and let $\Theta_n$ be the free module generated by Trees
with leaves labeled by $\n$.  Let $\Gamma_n$ be the free module generated by Graphs
with vertices labeled by $\n$.  Extend the configuration
pairing $\la \cdot \ra$ to one between $\Theta_n$
and $\Gamma_n$ by linearity.
\end{definition}

We next show that this pairing factors through canonical quotients of 
$\Theta_n$ and $\Gamma_n$ by Jacobi and Arnold identities.  Recall that Trees coincide
with elements of free non-associative algebras. 
For example, the first tree from \refF{basic}
is identified with $[[x_2, x_1], x_3]$.   In this language, if we replace a tree $T$ with a
bracket expression $B$, our pairing can be defined by the map $\beta_{G,B}$
which sends the edge between $i$ and $j$ to the innermost pair of brackets
which contains $x_{i}$ and $x_{j}$.

While it is traditional to define the Jacobi identity in the language of brackets,
namely that $[[A,B], C] + [[B,C],A] + [[C, A], B] = 0$ for any expressions $A$,$B$ and $C$,
we use the language of trees as follows.

\begin{definition}\label{D:jacobi}
\begin{enumerate}
\item A subtree of a Tree consists of a vertex  
and all edges and vertices
whose shortest path to the root goes through that vertex (that is, all edges and vertices
over that vertex).
\item A fusion of a Tree $T$ with another $S$ is the Tree obtained
by identifying the root edge of $S$ with a leaf edge of $T$, embedding $S$  
through a standard diffeomorphism of the upper-half plane with a
boundary-punctured disk disjoint from the rest $T$.
\item A Jacobi combination in $\Theta_n$ is a sum of three Trees
obtained by taking the tree $T$ from \refF{basic}, which has three leaves,
and fusing a tree $D$ to its root along with
three trees $A$, $B$ and $C$ to its leaves in three
cyclically-related orders.  See  \refF{jacobi}.
\item A symmetry combination is the sum of two Trees which are
isomorphic as graphs and have the same planar ordering of input edges at
each internal vertex but exactly one.\label{symmpart}
\item Let $J_n \subset \Theta_n$ be the submodule generated by 
Jacobi combinations and symmetry combinations.
\end{enumerate}
\end{definition}

\begin{proposition}\label{P:vanJ}
The pairing $\la \beta, \alpha \ra$ vanishes whenever $\alpha \in J_n$.
\end{proposition}

\begin{proof}
Vanishing on the symmetry combinations in $J_n$  is
immediate by our sign convention.
%, since the edge path $p_T(e)$ whose nadir
%is the vertex whose ordering is reversed will have its orientation
%in the plane reversed, whereas all other orientations will be the same
%for both Trees in the combination.

To check vanishing on a Jacobi combination, 
consider three Trees as in \refF{jacobi}.   In order for a Graph $\gamma$ 
to pair non-trivially with one of these trees (otherwise, the relation holds vacuously) 
the vertices $v_i$ and $w_i$ must be nadirs of edge paths.
There must then be precisely two distinct edge paths which begin at a leaf
of one of $A$, $B$ and $C$ and end in another.  Without loss of generality, we may 
assume that both of these edge paths begin or end at leaves in $A$.  

When paired with $\gamma$ the third term
of this sum is zero since the two edge paths share a common nadir, namely
$w_3$.  In the first
and second terms, the edge path between $A$ and $C$ has reversed its planar orientation.  
On the other hand, the edge path between $A$ and $B$ 
carries the same orientation in these two trees, as do all other edge paths of $\gamma$, 
which are either internal to $A$, $B$, $C$ or $D$ or which have one end in $D$.
Thus the signs of $\gamma$ paired with the first two trees will be opposite,
so that the sum of the three pairings  is zero.
\end{proof}

\begin{center}
\begin{minipage}{10cm}
$$\includegraphics[width=10cm]{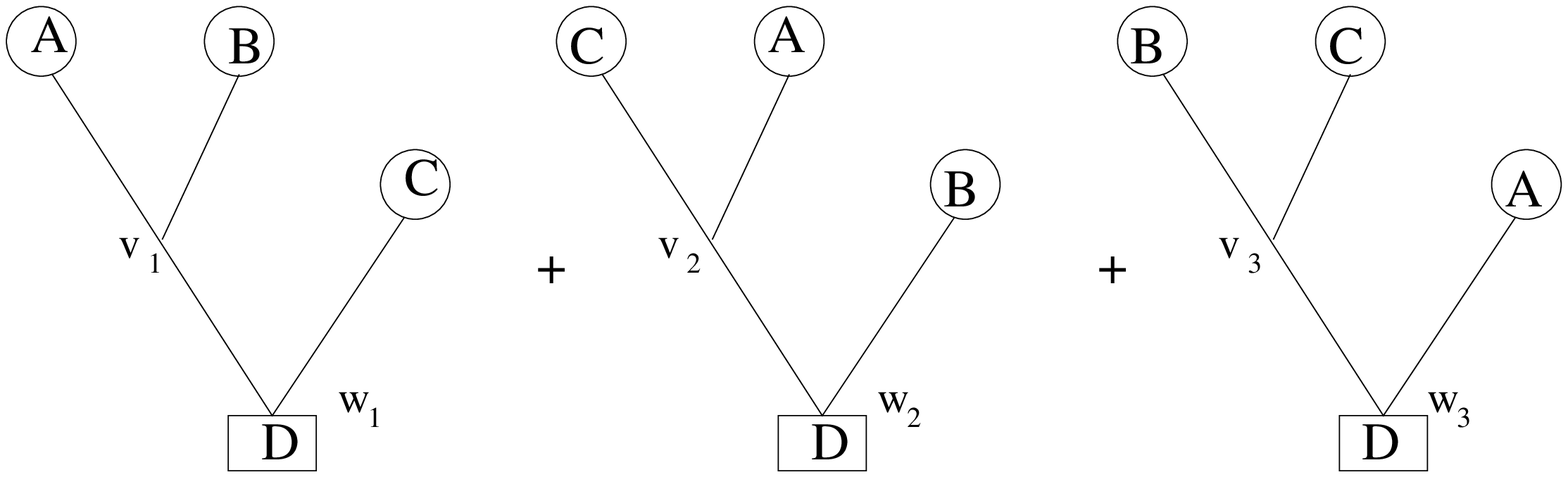}$$
\begin{mydiagram}\label{F:jacobi}
A Jacobi combination of Trees.
\end{mydiagram}
\end{minipage}
\end{center}
\medskip

We next consider a relation on Graphs which we will see as dual 
to the Jacobi identity.  

\begin{definition}
\begin{enumerate}
\item An Arnold combination in $\Gamma_n$
is the sum of three Graphs which differ only on the 
subgraphs pictured in \refF{dualjacobi}.  
\item A symmetry combination of Graphs
is the sum of two graphs which differ in the orientation of
exactly one edge.  
\item Let $I_n$ be the submodule of
$\Gamma_n$ generated by Arnold and symmetry combinations,
as well as Graphs with two edges which have the same vertices.
\end{enumerate}
\end{definition}

\begin{center}
\begin{minipage}{10cm}
$$\includegraphics[width=10cm]{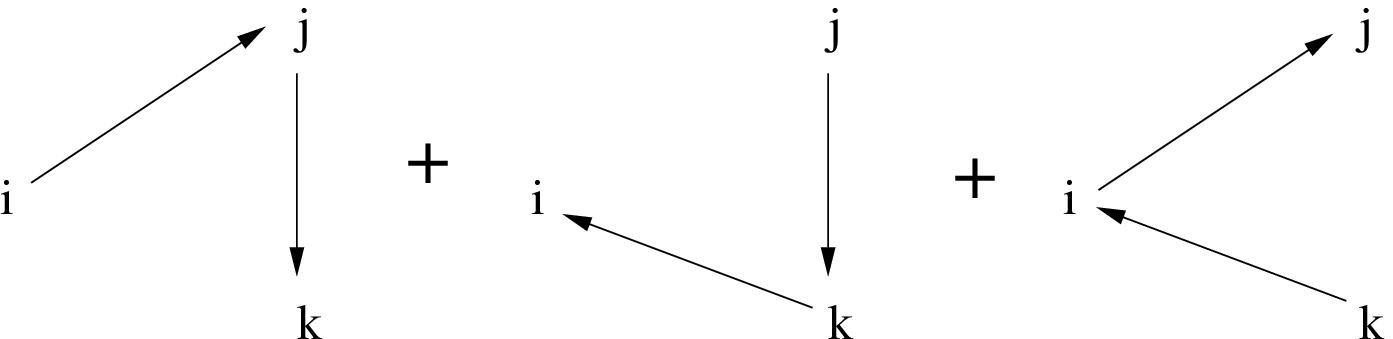}$$
\begin{mydiagram}\label{F:dualjacobi}
An Arnold combination of Graphs.
\end{mydiagram}
\end{minipage}
\end{center}
\medskip

Arnold combinations of graphs first occured, to our knowledge, in 
the computation of the cohomology of braid groups \cite{Arno68}.
They also appear in various forms of graph homology,
for example in Vassiliev's original work on knot theory  \cite{Vass92}.
We will see the Arnold combination of Graphs  as dual to the Jacobi 
combination of Trees, and thus in some contexts use the term Jacobi
combinations to refer to both.

\begin{proposition}\label{P:vanI}
The pairing $\la \beta, \alpha \ra$ vanishes whenever $\beta \in I_n$.
\end{proposition}

\begin{proof}
First, $\la G, T \ra$ vanishes whenever $G$ has two edges with the same
vertices since in this case $\beta_{G,T}$ cannot be a bijection. 
From our definition of $\tau$ it is immediate that the pairing of a 
tree $T$ with a symmetry combination of Graphs is zero.

To show that the pairing of an Arnold combination with 
a Tree $T$ vanishes, consider the nadir $v_{ij}$ of an 
edge path we now call $p_{ij}$ between leaves $i$ and $j$,
as well as the nadirs $v_{jk}$ and $v_{ki}$.  Two of these three nadirs must agree.  
Without loss of generality say  $v_{ij}$ and $v_{ki}$
agree, in which case the pairing
of $T$ with the third Graph in \refF{dualjacobi} is zero. 
The pairings with the first two Graphs in 
\refF{dualjacobi} differ by multiplication by $-1$, since the $p_{jk}$ 
appears with the same orientation in both cases but
$p_{ij}$ and $p_{ki}$ have the same nadir but different orientations.
The sum of these pairings is thus zero.
\end{proof}

\begin{definition}
Let $\Ln = \Theta_n/J_n$, and let $\Dn = \Gamma_n/I_n$.
\end{definition}

Propositions~\ref{P:vanJ} and~\ref{P:vanI} imply that the pairing
$\la, \ra$ passes to a pairing between $\Ln$ and $\Dn$, which
by abuse we give the same name.  %As is customary we use the term
%Jacobi and Arnold identities for equalities arising from setting a Jacobi or Arnold combination
%to zero, and the term anti-symmetry identity for equalities arising from
%symmetry combinations.
In $\Dn$ any Graph which has a cycle is zero, since
we may use the Arnold identity reduce to linear combinations
of Graphs with shorter cycles, ultimately until there are two edges 
which share the same vertices.
The modules $\Ln$ occur in many contexts, and in particular
are the entries of the $\Lie$ operad.

\begin{theorem}\label{T:perfect}
The pairing $\la, \ra$ between $\Ln$ and $\Dn$ is perfect.
\end{theorem}

Our proof uses reductions of these modules to particular bases.

\begin{definition}
\begin{enumerate}
\item The tall generators of $\Ln$ are represented by Trees for which the right
branch of any vertex is a leaf, and the leaf labeled by $1$ is leftmost.
\item The long generators of $\Dn$ are represented
by linear Graphs (all vertices but two are endpoints of exactly two edges)
with aligned orientations, with vertex $1$ as the initial endpoint.
\item Let $i_k$ be the label of the $k$th vertex from the left in a tall Tree or, respectively, in a 
long Graph, as in \refF{longntall}.  Given a permutation $\sigma \in \Sigma_n$ such
that $\sigma(1) = 1$ let $tT_\sigma$ and $lG_\sigma$ denote the tall Tree and long Graph,
respectively, with $i_k = \sigma(k)$.
\end{enumerate}
\end{definition}

%There are $(n-1)!$ of both of the tall and long generators.
%\begin{figure}
%\end{figure}

\begin{lemma}\label{L:longtallspan}
The tall generators span $\Ln$.  The long generators span $\Dn$.
\end{lemma}

\begin{proof}
Up to anti-symmetry identities, the tall generators of $\Ln$ are 
exactly those for which the leaf $1$ has the maximum height possible,
namely $n-1$.
We may reduce to such Trees by inductively applying the 
Jacobi, at each step getting a sum of two Trees each of which
has height of the leaf $1$ increased by one. 

In a similar spirit, we may start with any Graph generating $\Dn$
and first reduce to a linear combination of Graphs all of which have
a single edge from $1$  to another vertex and no other edges with vertex
$1$.  Indeed, if there is an edge between $1$ and $i$ and between
$1$ and $j$, we may after a change of orientation
use the Arnold identity to express it as a linear
combination of Graphs each with an edge instead between $i$ and $j$. 
This reduction decreases the number of edges with endpoint $1$.  
Given a graph with a single vertex from $1$ to
another vertex $i_2$, we repeat the above procedure  with $i_2$ in the place of $1$
to reduce until there is only
one edge from $i_2$ to some vertex $i_3$,
in addition to the one from $1$ to $i_2$.  
Repeating the procedure inductively reduces to linear graphs
with $1$ as an endpoint, and up to a sign we may change orientations
to align them.
\end{proof}

\begin{center}
\begin{minipage}{10cm}
$$\includegraphics[width=10cm]{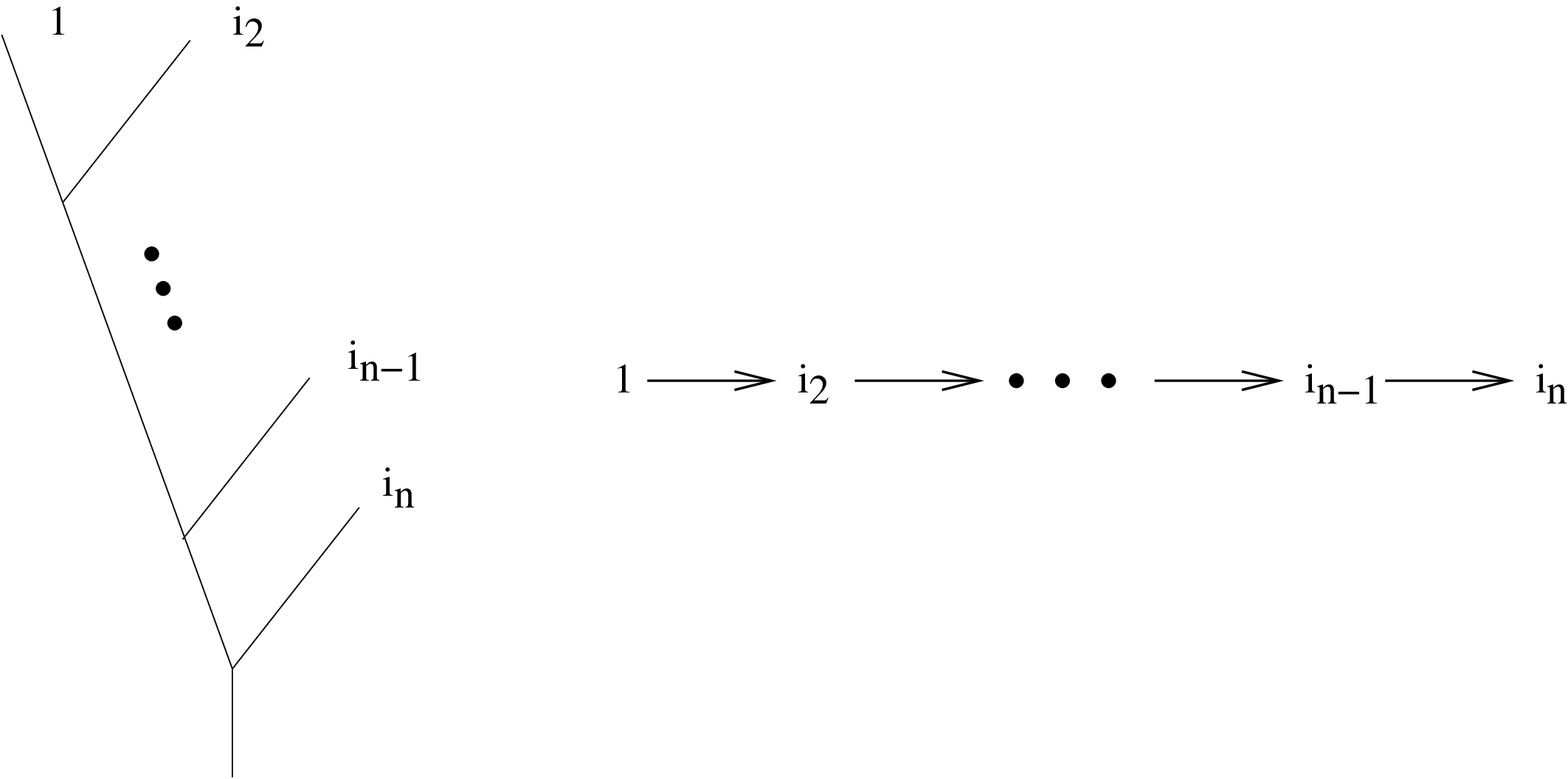}$$
\begin{mydiagram}\label{F:longntall}
A tall Tree and a long Graph.
\end{mydiagram}
\end{minipage}
\end{center}
\medskip

\begin{proof}[Proof of \refT{perfect}]
By the previous lemma, it suffices to show that the pairing is
perfect on the tall and long generators of $\Ln$ and $\Dn$ respectively.
By direct computation,
$\langle lG_\sigma, tT_\tau \rangle$ is one if $\sigma = \tau$
and zero otherwise, a perfect pairing.
\end{proof}
 
 In light of \refT{perfect}, we may view Jacobi, Arnold and anti-symmetry identities
 as arising as precisely the kernel of our pairing between Graphs
 and Trees.   

We may also deduce and extend Witt's classic calculation of a basis
for $\Ln$, without reference to the embedding of the free Lie algebra in the
corresponding free associative algebra.
  
 \begin{corollary}
 The tall Trees and long Graphs form bases for $\Ln$ and 
 $\Dn$ respectively, which are free of rank $(n-1)!$.   
 \end{corollary}
 
Moreover, we have a reduction method to these bases.

\begin{corollary}
Given $\alpha \in \Lie(n)$, $\alpha = \sum_{\substack{\sigma \in \Sigma_n \\ \sigma(1) = 1}}
   \la lG_\sigma, \alpha \ra tT_\sigma.$
 Similarly, given $\beta \in \Eil(n)$, $\beta = \sum_{\substack{ \sigma \in \Sigma_n \\\sigma(1) = 1}}
   \la \beta, tT_\sigma \ra lG_\sigma.$
 \end{corollary}
 
 For example, by looking at the tree $T$ which represents $\alpha = [[x_2, x_3],[x_1, x_4]]$, we 
 see that any long Graph which pairs with it non-trivially must start out with an edge
 from one to four, so the two possibilities are $1 \to 4 \to 3 \to 2$ and $1 \to 4 \to 2 \to 3$.
 By computing these pairings we have $\alpha = [[[x_1, x_4], x_3], x_2] - [[[x_1, x_4], x_2], x_3]$.
 
 \medskip 
 
 We reiterate that this duality between Lie trees and graphs modulo Arnold
 identities was first noticed by Tourtchine.  Using a recursive definition, it is developed in 
 \cite{Tour98}.  It is applied in a form close to ours in \cite{Tour04}.  Indeed, in that paper
 Section~2 gives various bases for the modules $\Dn$ (which are denoted $T^{+}_M$)
 and Section~5, in particular Statement~5.4, is devoted to the duality between $\Ln$
  (denoted $B^+_M$) and $\Dn$.

\section{The (even) Poisson configuration pairing}\label{S:even}
 
We next make the straightforward  passage 
to disconnected trees and graphs, which pertain
to Poisson algebras and configuration spaces.

 \begin{definition}\label{D:forest}
 \begin{enumerate}
\item Let $\Phi_n$ be the free module spanned by  unordered collections of Trees,
which we call Forests, with %(possibly fewer than $n$) 
leaves labeled by  $\n$. 
\item To a Forest $F$ associate a partition $\rho(F)$ of $\n$
by setting $i \sim j$ if $i$ and $j$ are leaves in the same tree.  Let 
$\Phi_n^P$ be the submodule spanned by all $F$ with $\rho(F) = P$.
We have $\Phi_n^P \cong \bigotimes_{S_i \in P} \Theta_{\# S_i}$.
\item Let $\Phi_n^k$ be the submodule spanned by all $F$ with a total of $k$
internal vertices.  We have that 
$$\Phi_n = \bigoplus_k \Phi_n^k = \bigoplus_P \Phi_n^P, \;\;\;
{\rm with} \;\;\; \Phi_n^k = \bigoplus_{P \; | \; \Sigma (\#S_i - 1) = k} \Phi_n^P.$$
\end{enumerate}
 \end{definition}

Thus $\Phi_n^{n-1}$ is isomorphic to $\Theta_n$, as is $\Phi_n^P$ where $P$
is the trivial partition.  

 \begin{definition}
 \begin{enumerate}
\item Let $\Delta_{n}$ be the free module spanned by unordered collections
of Graphs, which we call Diagrams,
 with vertices collectively labeled by $\n$.  Equivalently,
$\Delta_n$ is the free module spanned by  possibly disconnected oriented graphs.  
\item If $D$ is a Diagram, let $\rho(D)$ be the partition of $\n$ according
to the connected components of $D$.  Let $\Delta_n^P$ be the submodule
of $\Delta_n$ spanned by all Diagrams $D$ with $\rho(D) = P$.  We have that
$\Delta_n^P \cong \bigotimes_{S_i \in P} \Gamma_{\#S_i}$.
\item Let $\Delta_n^k$ be the submodule spanned by all $D$ with a total of $k$ edges.
We have that $$\Delta_n = \bigoplus_k \Delta_n^k = \bigoplus_P \Delta_n^P.$$
\end{enumerate}
 \end{definition}

 \begin{definition}\label{D:poispair}
 \begin{enumerate}
\item Extend the pairing $\la, \ra$ to $\Phi_n$ and $\Delta_n$ by setting
$\la D, F \ra$ to be zero unless $\rho(D) = \rho(F)$ and  pairing
$\Phi_n^P \cong \bigotimes_{S_i \in P} \Theta_{\# S_i}$ with 
$\Delta_n^P \cong \bigotimes_{S_i \in P} \Gamma_{\#S_i}$ through the
tensor product of the pairings between $\Theta_{\# S_i}$ and
$\Gamma_{\#S_i}$.
\item Define a Jacobi combination of Forests to be the sum of three forests whose 
component Trees are identical but for one component Tree of each, 
which together constitute a Jacobi combination of Trees.
Extend all definitions of Jacobi, Arnold, and symmetry combinations from Trees and Graphs
to Forests and Diagrams in similar fashion.
By abuse of notation, let $J_n$ be the submodule of $\Phi_n$  generated by Jacobi
and symmetry combinations, and let $I_n$ be the submodule of $\Delta_n$
generated by Arnold and symmetry combinations,
as well as Diagrams with two edges which have the same vertices.
\item Let $\Pn = \Phi_n/J_n$, and let $\DPn = \Delta_n/I_n$.   The submodules
$J_n$ and $I_n$ are generated by homogeneous elements, so these quotients
decompose as $$\Pn = \bigoplus_k \Pnk = \bigoplus_P \PnP \;\;\; {\text and} \;\;\;
\DPn = \bigoplus_k \DPnk = \bigoplus_P \DPnP,$$
with $\PnP \cong \bigotimes_{S_i \in P} \Lie(\#S_i )$ and
with $\DPnP \cong \bigotimes_{S_i \in P} \Eil(\# S_i)$.
\end{enumerate}
 \end{definition}
 
 Thus for example $\poi_{n-1}(n) \cong \Ln$. 

 This definition of  $\poi(n)$ is isomorphic to the more customary one using bracket
 expressions.  
 
 \begin{proposition}\label{P:bradef}
 $\poi(n)$ is isomorphic to the free module generated 
 expressions with two multiplications $[,]$ and $\cdot$
 in the variables $x_{1}, \ldots, x_{n}$, using each variable once, quotiented
 by anti-symmetry and the Jacobi identity in $[,]$, by commutativity in $\cdot$
 and by the Leibniz rule that $[,]$ defines a derivation with respect to  $\cdot$.
 \end{proposition}
 
\begin{proof}
We may use the Leibniz
 rule to reduce to expressions where there are no $\cdot$ multiplications which
 appear inside any bracket.  The resulting  
 $\cdot$ products of pure bracket expressions have no relations defined
 by the Leibniz rule, and are naturally represented by forests.
 The anti-symmetry, Jacobi and commutativity relations then translate
 exactly between these products of brackets and forests.
 \end{proof}

 The definition of the configuration pairing extends naturally to all  bracket
 expressions.
 
 \begin{definition}
 If $D$ is a Diagram and $B$ is a bracket expression define  the map $\beta_{D, B}$
 when possible by sending the edge between $i$ and $j$ to the innermost pair of brackets
 which contain $x_{i}$ and $x_{j}$.  If either no such pair of brackets exist
 or when $x_{i}$ and $x_{j}$ are multiplied by $\cdot$, internal to any brackets, we say
 that $\beta_{D, B}$ is not defined.
 
 Define the pairing $\la D, B \ra$ as $0$ if $\beta_{D,B}$ is not a bijection between
 the set of edges of $D$ and the set of bracket pairs of $B$ (in particular, if it is
 not defined), or $(-1)^{k}$, where $k$ is the number of edges $i \to j$ for which the 
 corresponding $x_{j}$ is to the left of $x_{i}$
 \end{definition}
 
 It is immediate that this definition of the pairing agrees with that of
 \refD{poispair} on bracket expressions which correspond to Forests.  
 Moreover, just as respecting Jacobi, Arnold and anti-symmetries
 was intrinsic in the Lie setting, respecting the Leibniz rule is intrinsic in this setting.
 
 \begin{proposition}\label{P:leib}
   $\la D, B \ra = \la D, B' + B'' \ra$, where $B$ is a bracket expression, $B'$ is obtained
   from $B$ by substituting $Y \cdot [X, Z]$ for $[X, Y \cdot Z]$, and $B''$ is obtained
   by substituting $[X, Y] \cdot Z$, for
   some sub-expressions $X, Y, Z$.
 \end{proposition}
 
 \begin{proof}
 The map $\beta_{D, B}$ will be a bijection if and only if exactly one of $\beta_{D, B'}$
 or $\beta_{D, B''}$ is a bijection, in which case the signs of the pairing will also agree.
 \end{proof}
 
 Finally, from \refT{perfect} and the 
 decompositions of $\poi(n)$ and $\DP(n)$ into (sums of) tensor products of $\Lie(\#S_i)$ and
 $\Eil(\# S_i)$, we immediately have the following.
 
 \begin{theorem}
 The configuration pairing $\la, \ra$ between $\Phi_n$ and $\Delta_n$
 descends to a perfect pairing between $\Pnk$ and $\DPnk$.
 \end{theorem}

\section{The odd Lie and Poisson configuration pairings}

There is a closely related pairing between graphs and trees whose sign is 
determined not by orientation of edges but by ordering.  Recall \refD{firstdef} as 
we make the following.

\begin{definition}\label{D:eoGraph}
\begin{enumerate}
\item An eoGraph is a connected graph with ordered edges and labeled 
vertices.
\item By abuse define $\beta_{G,T}$ and $\la G, T \ra$ for $G$ an eoGraph
and $T$ a Tree as in \refD{firstdef}.
\item Order the internal vertices of a Tree $T$ from left to right in accordance
with its embedding in the upper half-plane, so that if  a vertex $v$ sits over the left branch
of another vertex $w$ then $v<w$, and if $v$ sits over the right branch of $w$ 
then $w<v$.  \label{order}
\item Let $\sigma_{G,T}$ be the sign of the 
permutation defined through $\beta_{G,T}$ on the orderings of the edges
of $G$ and the internal vertices of $T$.  Define the pairing 
$\la G, T \ra$ as $\sigma_{G,T} \la G, T \ra_{(2)}$.
\item Let $o\Gamma_n$  be the free module generated by eoGraphs
with $n$ vertices, and extend $\la, \ra$ to a pairing between $\Theta_n$
and $o\Gamma_n$ by linearity.
\end{enumerate}
\end{definition}

The Jacobi and symmetry combinations of trees on which 
this pairing will vanish are the ones for graded Lie algebras generated
in odd degrees. 

\begin{definition}\label{D:oddrelns}
\begin{enumerate}
\item An odd Jacobi combination is a linear combination of three
Trees as in \refD{jacobi} but with each Tree having a coefficient as indicated
in \refF{ojacobi}, where $\#T$ is the number of leaves in a Tree $T$.
\begin{center}
\begin{minipage}{10cm}
$$\includegraphics[width=10cm]{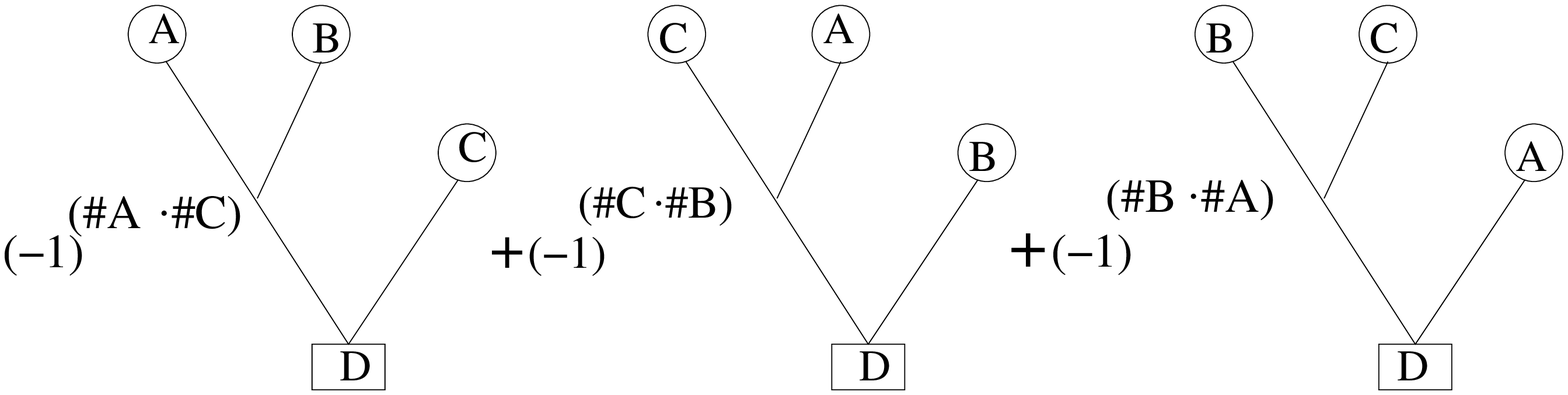}$$
\begin{mydiagram}\label{F:ojacobi}
An odd Jacobi combination of Trees.
\end{mydiagram}
\end{minipage}
\end{center}
\medskip

\item An odd symmetry combination  of two Trees is the sum of who which differ only
by switching left and right branches of a single vertex, as in part \ref{symmpart} of \refD{jacobi}, but
with a coefficient of $(-1)^{|A| |B|}$, where $A$ and $B$ are the subtrees emanating from
this vertex, for one of the resulting Trees. \label{I:oddsymm}
\item Let $oJ_n \subset \Theta_n$ be the submodule generated by 
odd Jacobi combinations and odd symmetry combinations.
\item An Arnold combination of eoGraphs
is the sum of three eoGraphs which differ only on the 
subgraphs pictured in \refF{odualjacobi}.  The two edges
in each of these subgraphs are ordered consecutively and
in the same position in all three eoGraphs; their ordering with
respect to each other is indicated by the labels of `I' and `II'.  

\begin{center}
\begin{minipage}{10cm}
$$\includegraphics[width=10cm]{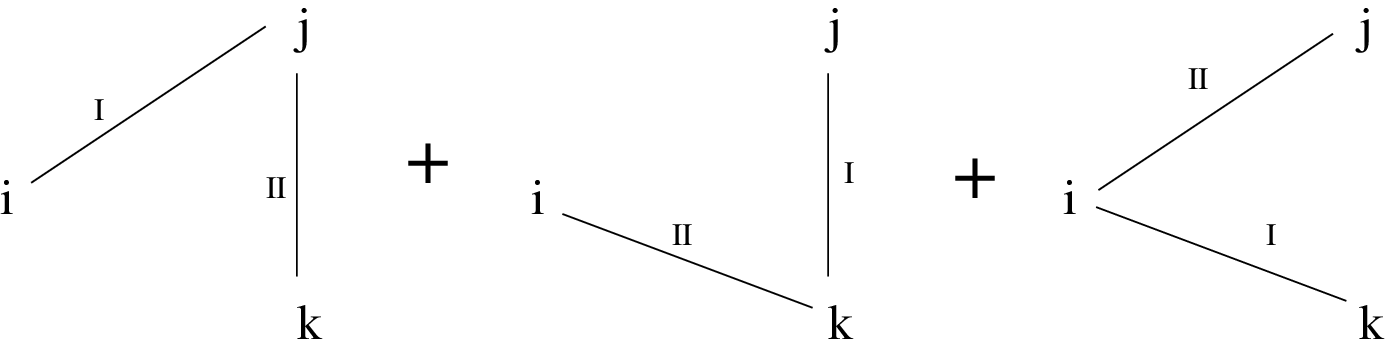}$$
\begin{mydiagram}\label{F:odualjacobi}
An Arnold combination of eoGraphs.
\end{mydiagram}
\end{minipage}
\end{center}
\medskip

\item A symmetry combination of eoGraphs
is a linear combination of two graphs which differ only in the
ordering of their edges, where one of the two eoGraphs has
a coefficient given by the sign of the permutation relating
these orderings.
\item Let $oI_n$ be the submodule of
$o\Gamma_n$ generated by Arnold and symmetry combinations,
as well as Graphs with two edges which have the same vertices.
\item Let $\olie(n)$ be $\Theta_n/oJ_n$, and let $\oeil(n)$ be 
$o\Gamma_n/oI_n$.
\end{enumerate}
\end{definition}

\begin{theorem}\label{T:oddLie}
The configuration pairing $\la, \ra$ passes to a perfect pairing 
between $\olie(n)$ and $\oeil(n)$.
\end{theorem}

\begin{proof}
Once we have shown that the configuration pairing vanishes
on $oJ_n$ and $oI_n$, the arguments from the proof of 
\refT{perfect} apply to show the pairing is perfect.    Indeed, the
tall generators of $\olie(n)$ share their definition with those of $\Lie(n)$,
and  we define the long generators of $o\Gamma_n$ to be the linear graphs
with vertex $1$ as an endpoint where the edges are ordered according
to their linear position.  There are still $(n-1)!$ tall generators
of $\olie(n)$ and long generators of $\oeil(n)$, which by inspection
pair perfectly.  The process of
\refL{longtallspan} applies almost verbatim, with only the coefficients
changing in the reduction process.

We show vanishing on odd Jacobi and symmetry combinations of Trees
and eoGraphs by straightforward computation.  
In pairing an eoGraph $\gamma$ with the two Trees
in an odd symmetry combination, only the signs of the two pairings
might differ since the Trees are isomorphic as graphs.  
The order of edges of $\gamma$ is fixed, so we 
consider the order of internal vertices of the two Trees, which
differ only around the vertex $v$ whose 
branches are the subtrees $A$ and $B$ of \refD{oddrelns}, part~\ref{I:oddsymm}.  
The transposition of the vertices in $A$ and $B$ is a composite of 
$(|A| - 1)|B| + |B|-1 = |A||B|-1$ transpositions.  When the resulting sign of $(-1)^{|A||B|-1}$
accounting for the difference between the two pairings is multiplied
by the $(-1)^{|A||B|}$ of \refD{oddrelns}  part (\ref{I:oddsymm}), we see
these two terms differ by their sign and thus cancel.

In pairing an eoGraph $\gamma$ with the three Trees 
in an odd Jacobi combination as in \refF{ojacobi},
we may as in the proof of \refP{vanJ} assume that there are 
precisely two distinct edges which begin at a leaf
of one of $A$, $B$ and $C$ and end in another, and that they both end at
leaves in $A$.  As before, the pairing of the third term with $\gamma$
is zero.  Because the ordering of edges in $\gamma$
is fixed, the difference in sign between the pairings
of $\gamma$ with the first and second trees is the sign of the
bijection between the  ordered sets of vertices in and between $A$, $B$
and $C$.  In the first Tree these vertices appear in the following order:
vertices in $A$, $v_1$, vertices in  $B$, $w_1$, vertices in $C$.
In the second Tree these vertices are in the order:
vertices in $C$, $v_2$, vertices
in $A$, $w_2$, vertices of $B$.  Under $\beta_{\gamma,T}$, $v_1$
and $w_2$ will correspond to the same vertex in $\gamma$,
as will $v_2$ and $w_1$.
Because there are $|T| - 1$ vertices in a subtree $T$,
the sign of this permutation is $-1$ to the power
$(|C| - 1)((|A| - 1) + 1 + (|B| -1) + 1) +
(|B| - 1) + 1 + (|A| - 1)$, which is equal to $(-1)^{|C|(|A| + |B|) - 1}$.
When, as in the definition of Jacobi combination, 
the first tree is multiplied by $(-1)^{|A| |C|}$ and the 
second by $(-1)^{|B| |C|}$ these two pairings will have
opposite signs and thus cancel.

The vanishing of the configuration pairing on odd symmetry combinations
of eoGraphs is immediate.  
To show that the pairing of an odd Arnold combination of graphs with 
a Tree $T$ vanishes, we may as in \refP{vanI} let $v_{ij}$ denote
the nadir of the edge path between leaves $i$ and $j$
and without loss of generality assume that $v_{ij}$ and $v_{ki}$
agree, in which case the pairing
of $T$ with the third Graph in \refF{odualjacobi} is zero. 
The pairings with the first two Graphs in 
\refF{dualjacobi} differ by multiplication by $-1$, since
$v_{ij} = v_{ki}$ is matched with edge I in the first eoGraph and
edge II in the second while $v_{jk}$ has the opposite matchings,
giving a sum of zero.
\end{proof}

We extend to the Poisson setting as in Section~\ref{S:even}. 
Recall \refD{forest} of the module of Forests.  There are completely
straightforward generalizations of these definitions where pre- or super-scripts
of $o$ decorate all of the modules.  We omit the repeated definitions.
Moreover, we may also use the language of expressions in variables
$x_{1}, \ldots, x_{n}$ using brackets and products, in which case the 
appropriately defined pairing respects the odd-graded Leibniz rule.

 From \refT{oddLie} and the definition of $\la, \ra$ on $\oPn$ and $\oDPn$ through
 their decompositions into (sums of) tensor products of $\olie(\#S_i - 1)$ and
 $\oeil(\# S_i - 1)$, we immediately have the following.
 
 \begin{theorem}
 The configuration pairing $\la, \ra$ between $\Phi_n$ and $o\Delta_n$
 descends to a perfect pairing between 
 $\oPnk$ and $\oDPnk$.
 \end{theorem}

\section{Functionals on free Lie algebras}\label{S:freelie}

 Working over a field in this section, let $V$ be a vector space  and $V^*$ be its dual.   
 Let $L(V)$ (respectively $L^o(V)$) be the (respectively oddly graded)
 free Lie algebra  generated by $V$, and $L_n(V)$ (respectively $L_n^o(V)$) be the 
 $n$th graded summand.  We construct functionals
 on $L(V)$ and $L^o(V)$ using elements of $V^*$.  
 First to set notation we recall that Trees
 may be used to define elements of $L(V)$ and $L^o(V)$.
 
 \begin{definition}
 Given $v_1, \ldots, v_n \in V$ and $T \in \Theta_n$, define $\mu_T(v_1, \ldots, v_n)$
 to be the product of $v_i$ according to $T$.  That is, translate $T$ into a bracket of
 free variables $x_i$ (as before \refD{jacobi}) substitute $v_i$ for $x_i$, and take the
 resulting product in $L(V)$ or $L^o(V)$.
 \end{definition}
 
 The $\mu_T(v_1, \ldots, v_n)$ span $L^o_n(V)$ (respectively $L^o(V)$).
 
 \begin{definition}\label{D:EV}
 Let $E(V)$, respectively $E^o(V)$, be the module spanned by Graphs
 (respectively eoGraphs), with vertices labeled by elements
 of $V$, up to linearity in each vertex and Arnold and symmetry relations.
 Given $w_1, \ldots, w_n \in V^*$ and $G \in \Gamma_n$ let  
 $\gamma_G(w_1, \ldots, w_n)$ denote the element of $E(V)$ 
 (respectively $E^o(V)$) where
 $w_i$ labels the $i$th vertex of $G$.
 \end{definition}
 
 \begin{definition}\label{D:liefnal}
 Define pairings between $L(V)$ and $E(V^*)$  (respectively
 $L^o(V)$ and $E^o(V^*)$) by 
 $$ \la \gamma_G(w_1, \ldots, w_n), \mu_T(v_1, \ldots, v_n) \ra
 = \sum_{ \sigma \in \Sigma_n} \left[ \la G, \sigma \cdot T \ra \prod_i w_{\sigma(i)} (v_i) \right],$$
 (and similarly in the odd case)
 where $\sigma \in \Sigma_n$ acts on $T$ by permuting the labels of its leaves.
 Extend to all of $L(V)$ linearly.
 \end{definition}
 
 \begin{proposition}
 The pairing of \refD{liefnal}  is well-defined.
 \end{proposition}
 
 \begin{proof}
Note that $\mu_T(v_1, \ldots, v_n) = \mu_{\tau \cdot T}(v_{\tau(1)}, \ldots, v_{\tau(n)})$, 
and similarly for  $\gamma_G$ for any permutation $\tau$.   We  compute immediately 
that the pairing gives
the same value on any of these representations.   With this  established, the only 
other equalities to check arise from Jacobi and (anti-)symmetry identities in $L(V)$
along with Arnold and symmetry identities in $E(V)$,
which with freedom to permute the leaves of $T$ and the $v_i$ all follow
from the vanishing of $\la, \ra$ on the corresponding identities of (o)Graphs and
Trees.
 \end{proof}
 
 Take for example the free Lie algebra on two letters, so that $V$ is spanned
 by say $a$ and $b$.  Let $T$ be the first tree in \refF{basic} and let
 $v_1 = a, v_2 = b, v_3 = b$, so that $\mu_T(v_1, v_2, v_3) = [[b,a]b]$.
 Let $G$ be as in \refF{basic}, and compute that  $\la \gamma_G(a^*, b^*, b^*), [[b,a]b] \ra$
 is one, since the only term in the sum of \refD{liefnal} which is non-zero is the
 one with $\sigma = id$.  All other terms vanish either because $\prod w_{\sigma(i)}(v_i)$
 will be zero, since some term such as $b^*(a)$ will occur, or in the case
 that $\sigma$ transposes $2$ and $3$ because the resulting $\la G, \sigma \cdot T \ra$
 is zero, as computed in \refF{basic}.
 
 Melancon and Reutenaur give a pairing equivalent to this one for graded
 Lie algebras generated in odd degrees in \cite{MeRe96}.  
 Their main result is the following.
 
 \begin{theorem}[Theorem~4.1 of \cite{MeRe96}]
 The pairing between $L^o(V)$ and $E^o(V^*)$ is perfect.
 \end{theorem}
 
 Thus $E(V^*)$ is a model for the free Lie coalgebra on $V^*$.
We extend this result to graded Lie algebras, without  using embeddings
of free Lie algebras in associative algebras, in \cite{SiWa06}.
 
 \section{Product and coproduct structures}
 
 The modules of Diagrams $\Delta_n$ and $o\Delta_n$ carry  a multiplication
 which is elementary.  
 
 \begin{definition}
 Given Diagrams (respectively oDiagrams)
 $D_1$ and $D_2$ their product  $D_1 \cdot D_2$ is the Diagram $D$ whose edges are the
 union of the edges in $D_1$ and $D_2$, carrying the same orientations (or 
 respectively having the edges of $D_1$ occur in their given order before 
 the edges of $D_2$ occur in their given order).
 \end{definition}
 
 Under this product $\Delta_n$ and $o\Delta_n$ are commutative
 (respectively graded commutative, with $o\Delta_n^k$ in degree $k$)
 rings with unit, namely the Diagram with no edges.
 
 \begin{proposition}
 The multiplication  on $\Delta_n$ and $o\Delta_n$ passes to  
 the quotients $\DPn$ and $\oDPn$.
 \end{proposition}
 
 Just as the product on $\DPn$ and $\oDPn$ are defined through
 unions of edges, there is a coproduct on $\Pn$ and $\oPn$ 
 defined through partitions of internal vertices.
 
 \begin{definition}
 A partition of a Forest $F$ with $n$ leaves and $k$ internal vertices
 is a pair of Forests $F_1, F_2$
 each with $n$ leaves and with $\ell_1$ and $\ell_2$ internal vertices
 where $\ell_1 + \ell_2 = k$, defined as follows.
 \begin{enumerate}
 \item Partition the set of internal vertices of $F$ into two sets,
 $S_1$ and $S_2$.
 \item For each vertex in $S_1$, choose one leaf
 above each of the two branches of that vertex. \label{step2}
 \item Take the smallest subgraph $G_1$ of $F$ containing all of the 
 vertices in $S_1$ along with all of the leaves chosen in the previous step.
 \item Obtain $F_1$ by replacing pairs of edges of $G_1$ connected by a bivalent 
 vertex with a single edge, adding a root edge to the lowest vertex of each
 connected component, and adding a single-edge Tree for each leaf 
 vertex not chosen in step~\ref{step2}.
 \item Repeat the previous three steps using $S_2$ instead of $S_1$ to obtain
 $F_2$.
 \end{enumerate}
 \end{definition}
  
 \begin{center}
\begin{minipage}{14cm}
$$\includegraphics[width=14cm]{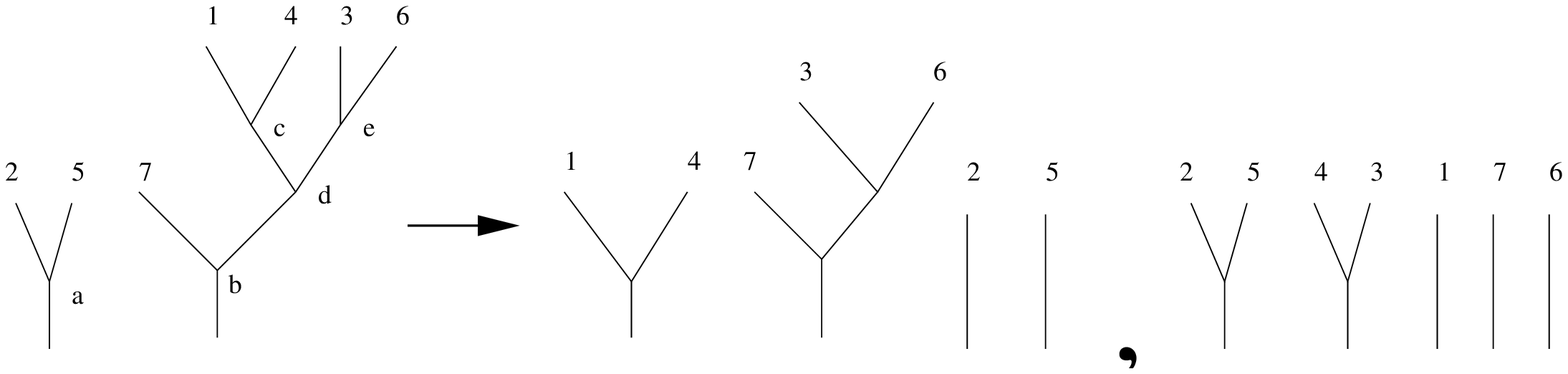}$$
\begin{mydiagram}\label{F:pruning}
A partition of a Forest $F$ into $F_1, F_2$.  The internal vertices chosen for $F_1$
are those labeled by $b, c, e$, and for $F_2$ are $a, d$.
\end{mydiagram}
\end{minipage}
\end{center}

\begin{definition}
Define a coassociative, cocommutative coproduct 
$c: \Phi^k_n \to \bigoplus_{\ell_1 + \ell_2 = k} \Phi^{\ell_1}_n \otimes \Phi^{\ell_2}_n$
by $c(F) = \sum_{(F_1, F_2) \in P_F} F_1 \otimes F_2$, where
$P_F$ is the set of partitions of $F$. 
\end{definition}
 
 \begin{theorem}\label{T:coprodpair}
 Let $F \in \Phi_n^k$ and let $G_1$ and $G_2$ in $\Delta_n^{\ell_1}$ and $\Delta_n^{\ell_2}$
 (or $o\Delta_n^{\ell_1}$ and $o\Delta_n^{\ell_2}$ respectively)
 with ${\ell_1} + {\ell_2} = k$.  Then $\la G_1 \cdot G_2, F \ra = 
 \la G_1 \otimes G_2, c(F)  \ra_*,$ where $\la, \ra_*$ is the direct sum of
 tensor products of configuration pairings between 
 $\bigoplus_{{\ell_1} + {\ell_2} = n} \Phi_n^{\ell_1} \otimes \Phi_n^{\ell_2}$ and
 $\bigoplus_{{\ell_1} + {\ell_2} = n} \Delta_n^{\ell_1} \otimes \Delta_n^{\ell_2}$ (or respectively 
 $\bigoplus_{{\ell_1} + {\ell_2} = n} o\Delta_n^{\ell_1} \otimes o\Delta_n^{\ell_2}$).
 \end{theorem}

 \begin{proof}
If $\beta_{G_1 \cdot G_2, F}$ is a bijection then 
 the nadirs of the paths $p_F(e)$ for $e$ in $G_1$ and the nadirs for $e \in G_2$
 partition the internal vertices of $F$ into two sets.   If we remember the
 vertex labels of the edges $e$ we also have
 a choice of two leaves over each internal vertex of $F$, which
 gives rise to a unique partition of $F$ into say $(\phi_1, \phi_2)$.  From the definition
 of partition of a Forest we see that  
 $\beta_{G_1, \phi_1}$ and $\beta_{G_2, \phi_2}$ are bijections.  
 Moreover, the sign $\tau_{G_1 \cdot G_2, F}$ equals the product 
 $\tau_{G_1, \phi_1} \cdot \tau_{G_2, \phi_2}$ (respectively $\sigma_{G_1 \cdot G_2, F}
 = \sigma_{G_1, \phi_1} \sigma_{G_2, \phi_2}$).  Finally, for no other
 partition of $F$ into some $(F_1, F_2)$ will 
 $\beta_{G_1, F_1}$ and $\beta_{G_2, F_2}$ be bijections.   In all
 such cases there will be some edge in either $G_1$ or $G_2$ 
 whose endpoints correspond to leaves which are in different
 connected components of $F_1$ or $F_2$.  
 Thus 
 \begin{multline*}
 \la G_1 \otimes G_2, c(F) \ra_* = 
 \la  G_1 \otimes G_2, \sum_{(F_1, F_2) \in P_F} F_1 \otimes F_2 \ra_* =
 \la G_1 \otimes G_2 , \phi_1 \otimes \phi_2 \ra_* = \\
 \la G_1, \phi_1 \ra \la G_2, \phi_2 \ra = \tau_{G_1, \phi_1} \cdot \tau_{G_2, \phi_2}
 = \tau_{G_1 \cdot G_2, F} = \la G_1 \cdot G_2, F \ra,
 \end{multline*}
 or similarly with $\sigma$'s replacing $\tau$'s in the odd setting.
 
 If $\beta_{G_1 \cdot G_2, F}$ is not a bijection, with say an internal vertex
 $v$ not in its image, then for any partition $(F_1, F_2)$ of $F$, $v$ cannot
 be in the image of $\beta_{G_1, F_1}$ or $\beta_{G_2, F_2}$, so both
 $\la G_1 \cdot G_2, F \ra$ and $\la G_1 \otimes G_2, c(F) \ra_*$ are zero.
  \end{proof}

We conclude this section with the following.
 
 \begin{corollary}\label{C:coprodvan}
 The coproduct $c$ passes from $\Phi_n$ to its quotients $\Pn$ and $\oPn$.  Thus the equality
 of   \refT{coprodpair} holds for $F \in \Pn$ and $G_i \in \DPn$ or $F \in \oPn$ and $G_i \in \oDPn$.
 \end{corollary}
 
 \begin{proof}[Proof of \refC{coprodvan}]
 Assume that $F$ is in the Jacobi and symmetry submodule $J_n$ (respectively $oJ_n$).
 By \refT{coprodpair}, $\la G_1 \otimes G_2,  c(F)  \ra_* = \la G_1 \cdot G_2, F  \ra = 0$,
 since the configuration pairing vanishes on $F$.  But the pairing $\la , \ra_*$ is perfect,
 so $c(F) = 0$.
 \end{proof}
  
 \section{Operad structures}\label{S:op}

It is well-known that the $\Ln$ and $\Pn$ assemble to form operads.  
In this section we determine the linearly dual cooperad structure on the $\Dn$ and $\DPn$.
Similar ideas were developed in \cite{KoSo00} for graph complexes.
For simplicity we restrict attention to the non-$\Sigma$ operad
structure, giving the following definitions in order to set notation.

\begin{definition}
\begin{enumerate}
\item Rooted planar trees, which we call rp-trees, share much of their definition
with Trees, but are unlabeled and not restricted to have only either trivalent or 
univalent vertices.
% and have leaves which are always
%labeled by $\n$ from left to right in the upper half plane.

\item Given an rp-tree $\tau$ and a set of edges $E$ the 
contraction of $\tau$ by $E$ is the rp-tree $\tau'$ obtained by, for 
each edge $e \in E$, identifying its two vertices 
(altering the embedding  only in a small neighborhood 
of $e$) and removing $e$ from the set of edges. 

\item Let $\Upsilon$ denote the category of rp-trees, in which there is
a unique morphism $f_{\tau, \tau'}$ from $\tau$ to $\tau'$,
if $\tau'$ is the
contraction of $\tau$ along some set of  non-leaf edges $E$.  
Let  $\Upsilon_n$ denote the full subcategory
of rp-trees with $n$ leaves. 

\item  Each $\Upsilon_n$ has a terminal object, namely the unique tree with one vertex,
called the $n$th corolla $\gamma_n$ as in \cite{MSS02}.  We allow for the tree
$\gamma_0$ which has no leaves, only a root vertex, and is
the only element of $\Upsilon_0$.

\item An edge is called
redundant if one of its vertices is bivalent.
 For a vertex $v$ let $|v|$ denote its valence minus one.  
 
\end{enumerate}
\end{definition}

\begin{definition}\label{D:op}
A non-$\Sigma$ operad is a functor $\Op$ from $\Upsilon$ to a symmetric
monoidal category $(\mC, \odot)$ which satisfies the following axioms.

\begin{enumerate}
\item $\Op(\tau) = \odot_{v \in \tau} \Op(\gamma_{|v|})$.  \label{1}
\item $\Op(\gamma_1) = \one_{\mC} = \Op(\gamma_0)$. \label{2}
\item If $e$ is a redundant edge and $v$ is its terminal vertex
then under the  decomposition of axiom (\ref{1}) $\Op(c_{\{e\}})$ is the identity map on
$\odot_{v' \neq v} F(\gamma_{v'})$ tensored  with the isomorphism
$(\one_{\mC} \odot -)$. \label{3}
\item If $\mu$ is a subtree of $\tau$, by which we mean a collection of
vertices and their branches which is itself a tree, and if $f_{\mu, \mu'}$ and $f_{\tau, \tau'}$
contract the same set of edges, then under the
decomposition of (\ref{1}), $F(f_{\tau,\tau'}) = F(f_{\mu, \mu'}) \odot id$. \label{4}
\end{enumerate}
\end{definition}

By axiom~(\ref{4}),
the values of $\Op$ on morphisms may be computed by
composing morphisms on subtrees, so we may identify some 
subset of basic morphisms through which all morphisms
factor.  The basic class we consider is 
that of all morphisms $\tau \to \gamma_n$ where $\gamma_n$ is a 
corolla. This class includes the $\circ_i$ operations
and May's structure maps.

\begin{definition}
The module $\Theta = \bigoplus_i \Theta_i$ forms an operad 
which associates to the morphism $\tau \to \gamma_n$ in $\Upsilon$ the homomorphism
$f_\tau$ sending $\bigotimes_{v_i} T_{v_i}$, where $v_i$ ranges over internal
vertices in $\tau$ and $T_{v_i} \in \Theta_{|v_i|}$, to the tree $S \in \Theta_n$,
called the grafting of the $T_{v_i}$ and obtained as a quotient of them  as follows.
If $v_i$ in $\tau$, is the other vertex of
the $k$th branch (in the planar ordering of edges) of $v_j$, we
identify the root edge of $T_{v_i}$ with the $k$th leaf
of $T_{v_j}$.  Label the vertices of $S$ by elements of $\n$
according to the total ordering where leaf $\ell$ is less than
leaf $m$ if they both sit over some $T_{v_i}$ and the leaf of $T_{v_i}$
over which $\ell$ sits has a smaller label than that over which $m$ sits.

This operad structure passes immediately
to the quotient $\Lie = \bigoplus_i \Lie(i)$, known
as the Lie operad.   In the odd setting of $\olie = \bigoplus_i \olie(i)$,
$f_\tau$ sends $\bigotimes_{v_i} T_{v_i}$ to $({\rm sign} \rho) S$,
where $S$ is as above and $\rho$ is the permutation which relates
the order of internal vertices as they occur in $\bigotimes_{v_i} T_{v_i}$
with the order of the corresponding internal vertices of $S$.
\end{definition}

Through the configuration pairing duality, we know that 
$\Eil = \bigoplus_i \Eil(i)$ form a cooperad, which we understand explicitly as follows.   

\begin{definition} \label{D:mB}
Label both the leaves of an rp-tree and the 
branches of each internal vertex $v$ with elements of
$\n$ and $\mathbf{|v|}$ respectively,
from left to right using the orientation in the upper half plane.  
To an rp-tree $\tau$ with $n$ leaves and
two distinct integers $j,k \in \n$ let $v$ be the nadir
of the shortest path between leaves labelled $i$ and $j$ and define
$J_v(j), J_v(k)$ to be the labels of the branches of $v$ over
which leaves $j$ and $k$ lie.

The module $\Gamma = \bigoplus_i \Gamma_i$,
forms a cooperad  which associates to the morphism $\tau \to \gamma_n$ 
the homomorphism $g_\tau$ sending $G \in \Gamma$ to 
$\bigotimes_{v_i} G_{v_i}$, as $v_i$ ranges over internal
vertices in $\tau$.  The graph $G_{v_i} \in \Gamma_{|v_i|}$, 
is defined by having for each edge in $G$, say between $j$ and $k$, an
an edge between $J_v(j)$ and $J_v(k)$ in $G_v$.
The module  $o\Gamma = \bigoplus_i o\Gamma_i$ similarly forms
a cooperad with structure map $g_\tau$ sending $G$
to $ {\rm{sign}} \pi \bigotimes_{v_i} G_{v_i},$ with $G_{v_i}$
as above, with an ordering of its edges given by the order
of the edges in $G$ which give rise to them, and $\pi$ is the
permutation relating this order on all of the edges in 
 $\bigotimes_{v_i} G_{v_i}$ to the ordering within $G$.
\end{definition}

For example, by a small abuse we may consider the second tree of the 
Forest $F$ from Figure~\ref{F:pruning}.  If $j, k = 3, 4$ then
$v$ is the vertex labeled $d$, and $J_{d}(3) = 2$ while $J_{d}(4) = 1$.

This definition is closely related to 
the choose-two operad, introduced in Section~2.2 of \cite{Sinh04}.
This operad structure is illustrated
in the more general setting of Forests and Diagrams in 
Figure~\ref{F:operad}.

If we consider $E(V)$ as in \refD{EV}, then the canonical cooperad action
induces a Lie coalgebra structure under which  a tree $G$ maps to 
$\sum_{e \in G} G' \otimes G'' - G'' \otimes G'$, 
where $G'$ and $G''$ are the sub-trees obtained by removing $e$.
We use this Lie coalgebra structure in \cite{SiWa06}.

\begin{theorem}\label{T:operadpair}
Let $\tau$ be an rp-tree and let $T_{v_i} \in \Theta_{|v_i|}$, with
$v_i$ ranging over the internal vertices of $\tau$.  Let $G \in \Gamma_n$ or
$o\Gamma_n$, where $n = \sum |v_i|$.  Then
$\la G, f_\tau(\bigotimes T_{v_i}) \ra = \la g_\tau(G), \bigotimes T_{v_i} \ra_\otimes$,
where $\la , \ra_\otimes$ denotes the tensor product of (respectively
even or odd) configuration pairings.
\end{theorem}

\begin{proof}
A vertex $w$ in the $T_v$ subtree of 
$T = f _\tau(\bigotimes T_{v_i})$ is in the image of $\beta_{G,T}$
if and only if there is an edge in $G$ with vertices whose labels coincide with those of 
one leaf above
the left branch of $w$ and one leaf above the right branch.  Such leaves
correspond to the 
leaves of $\tau$ which lie above the corresponding edges of $v$. 
We see that this is also the condition for  
$w$ to be in the image of $\beta_{G_v, T_v}$ as well.  

To determine the signs in the even setting note that
in pairing with both $T$ and $\bigotimes T_{v_i}$, the orientation
of an edge path as it passes through $w$ depends only on 
the leaves of $T_v$ connected by that edge, which in turn in both
cases only depends
on the leaves of $\tau$ connected by that edge path.   When $G \in o\Gamma_n$,
the signs of the permutations $\rho$ and $\pi$
in the definition of $f_\tau$ and $g_\tau$ relate the signs
of $\la G, f_\tau(\bigotimes T_{v_i}) \ra$ and 
$\la g_\tau(G), \bigotimes T_{v_i} \ra_\otimes$.

In summary $\beta_{G,T}$
is a bijection if and only if all of the $\beta_{G_v, T_v}$ are, and signs
agree, establishing the result.   
\end{proof}

The argument of \refC{coprodvan} adapts to this setting to give the following main result.

 \begin{corollary}\label{C:operadvan}
 The (co)operad structure map $g_\tau$ passes from 
 $\Gamma_n$ and $o\Gamma_n$ to their quotients $\Dn$ and
 $\Eil^o(n)$ respectively.  
 %Thus the equality
 %of   \refT{operadpair} holds for $T_{v_i}  \in \Lie(|v_i|)$ and $G \in \Dn$.
 \end{corollary}

Finally, we treat the case of Forests and the Poisson operad.
Here we use \refP{bradef}, choosing to describe the operad structure
in terms of bracket expressions.

\begin{definition}
\begin{enumerate}
\item  Let $f: \tau \to \gamma_{n}$ be a morphism in $\Upsilon$ in which $\tau$
is a tree with one internal vertex over the $i$th root edge.  
Such morphisms give rise to what are known
as $\circ_{i}$-operations, which generate an operad structure.
Define an operad structure on $\poi = \oplus \poi(n)$ by sending $f$
to the map $\poi(n) \otimes \poi(m) \to \poi(n+m -1)$ where $B_{1} \otimes B_{2}$
is sent to the bracket expression in which the variables in $B_{2}$ are re-labeled
from $x_{i}$ to $x_{m}$, the variables $x_{j}$ in $B_{1}$ with $j>i$ are re-labeled
by $x_{j+m-1}$, and then $B_{2}$ is substituted for $x_{i}$
in $B_{1}$.

\item Apply \refD{mB} of a (co)operad structure on the Graph module $\Gamma$
verbatim to define such a structure on the Diagram module $\Delta = \bigoplus_i \Delta_i$.
\end{enumerate}
\end{definition}

\begin{theorem}\label{T:operadpair2}
Let $\tau$ be an rp-tree and $f_{\tau}$ the morphism from $\tau$ to the corresponding
corolla.  Let $B_{v_i}$ be bracket expressions in $|v_{i}|$ variables, where
$v_i$ ranges over the internal vertices of $\tau$.  Let $D \in \Delta_n$ or
$o\Delta_n$,
where $n = \sum |v_i|$.  Then
$\la D, f_\tau(\bigotimes B_{v_i}) \ra = \la g_\tau(D), \bigotimes B_{v_i} \ra_\otimes$,
where $\la , \ra_\otimes$ denotes the tensor product of (respectively
even or odd) configuration pairings.
\end{theorem}

The proof is entirely analogous to that of \refT{operadpair}.   The proof would be more
involved if we had used Forests instead of bracket expressions to define the operad
structure, because the Leibniz rule must then be used to compute the final result of
a structure map.  But by \refP{leib}, the configuration pairing works perfectly 
well for bracket expressions, for which the proof is straightforward.

 \begin{corollary}\label{C:operadvan2}
 The cooperad structure map $g_\tau$ passes from 
 $\Delta_n$ and $o\Delta_n$ to their quotients $\DPn$ and
 $\oDPn$ respectively.  
 \end{corollary}
 
 This cooperad structure is more manageable than the operad structure on 
 $\poi$, for which the Leibniz rule is needed to reduce to any basis.
 For example, in \cite{Sinh06} where we establish the classical result that the homology of
 the little disks operads are $\poi$ or $\poi^{o}$, it is simpler to work with 
 cohomology and show that the cooperad structure there agrees $\DP$ or $\DP^{o}$.

 We end with a small illustration of \refT{operadpair2}.

\begin{center}
\begin{minipage}{15cm}
{\psfrag{tau}{$\tau$}
\psfrag{ftau}{$f_\tau$}
\psfrag{gtau}{$g_\tau$}
$$\includegraphics[width=12cm]{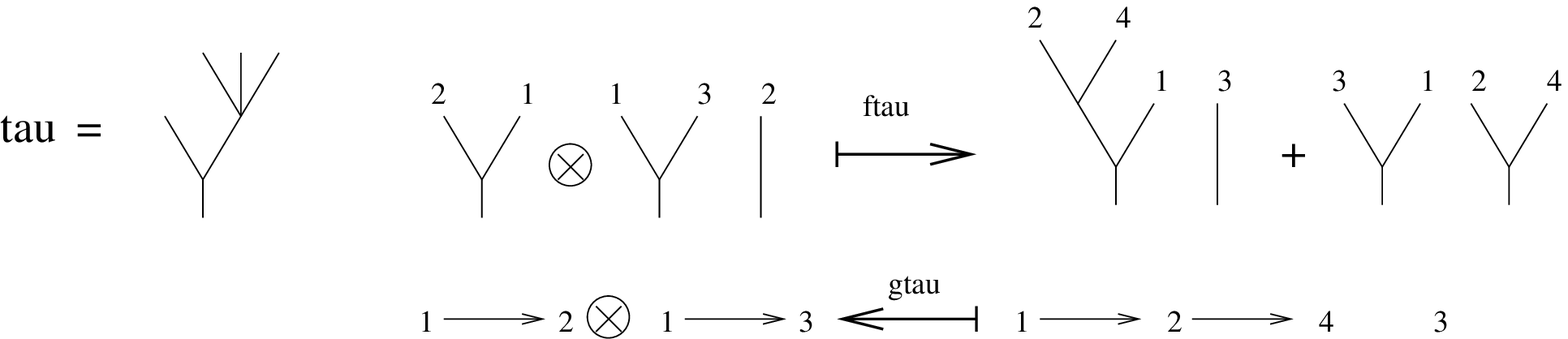}$$}
\begin{mydiagram}\label{F:operad}
Examples of the operad structure maps $f_\tau$ and $g_\tau$.  
Both $\la D, f_\tau(\bigotimes F_{v_i}) \ra$ and $\la g_\tau(D), \bigotimes F_{v_i} \ra_\otimes$
are equal to $-1$.
\end{mydiagram}
\end{minipage}
\end{center}

\end{document}